\documentclass{amsart}
\usepackage{graphicx}
\usepackage{amssymb,amscd,amsthm,amsxtra}
\usepackage{latexsym}
\usepackage{epsfig}

\newtheorem{thm}{Theorem}[section]
\newtheorem{cor}[thm]{Corollary}
\newtheorem{lem}[thm]{Lemma}

\theoremstyle{definition}

\theoremstyle{remark}
\newtheorem{rem}[thm]{Remark}
\numberwithin{equation}{section}

\newcommand{\R}{\mathbb R}

\newcommand{\eps}{\varepsilon}

\newcommand{\p}{\partial}

\newcommand{\comment}[1]{}

\begin{document}

\title[Global $W^{2,p}$ estimates]{Global $W^{2,p}$ estimates for the Monge-Ampere equation}
\author{O. Savin}
\address{Department of Mathematics, Columbia University, New York, NY 10027}
\email{\tt  savin@math.columbia.edu}

\thanks{The author was partially supported by NSF grant 0701037.}

\begin{abstract}
We obtain global $W^{2,p}$ estimates for the Monge-Ampere equation under natural assumptions on the domain and boundary data.
\end{abstract}
\maketitle

\section{Introduction}

Interior $W^{2,p}$ estimates for strictly convex solutions for the Monge-Ampere equation were obtained by Caffarelli in \cite{C} under the necessary assumption of small oscillation of the right hand side. The theorem can be stated as follows.

\begin{thm}[Interior $W^{2,p}$ estimates]
Let $u: \overline \Omega \to \R$,
\begin{equation}
u=0 \quad \mbox{on $\p \Omega$},
\end{equation}
be a continuous convex solution to the Monge-Ampere equation
\begin{equation}\label{MA}
\det D^2 u=f(x) \quad \mbox{in $\Omega$}, \quad \quad 0< \lambda \le f \le \Lambda,
\end{equation}
for some positive constants $\lambda$, $\Lambda$.

 For any $p$, $1<p<\infty$ there exists $\eps(p)>0$ depending on $p$ and $n$ such that if
 $$B_\rho \subset \Omega \subset B_{1/\rho}$$
 and
\begin{equation}\label{f_osc}
\sup_{|x-y|\le \rho}|\log f(x)-\log f(y)|\le \eps(p),
\end{equation}
for some small $\rho>0$, then $$\|u\|_{W^{2,p}(\{u<-\rho\})} \le C,$$
and $C$ depends on $\rho$, $\lambda$, $\Lambda$, $p$ and $n$.
\end{thm}

In this short paper we obtain the global $W^{2,p}$ estimates under natural assumptions on the domain and boundary data.

\begin{thm}[Global $W^{2,p}$ estimates]\label{w2p}
Let $\Omega$ be a convex bounded domain and $u:\overline \Omega \to \R$ be a Lipschitz continuous convex solution of the Monge-Ampere equation \eqref{MA}. Assume that $$u|_{\p \Omega}, \quad \p \Omega \in C^{1,1},$$ and there exists $\rho>0$ small such that $f$ satisfies \eqref{f_osc}. If $u$ separates quadratically on $\p \Omega$ from its tangent planes i.e
\begin{equation}\label{quad_sep}
u(y)-u(x)-\nabla u(x) \cdot (y-x) \ge \rho |x-y|^2 \quad \quad \forall x,y \in \p \Omega,
\end{equation}
then
$$\|u\|_{W^{2,p}(\Omega)} \le C,$$ with $C$ depending on $\|\p \Omega\|_{C^{1,1}}$, $\| u|_{\p \Omega}\|_ { C^{1,1}}$, $\|u\|_{C^{0,1}}$, $\rho$, $\lambda$, $\Lambda$, $p$ and $n$.
\end{thm}

\begin{rem}
The gradient $\nabla u(x)$ when $x \in \p \Omega$ is understood in the sense that $$y_{n+1}=u(x)+\nabla u(x)\cdot (y-x)$$ is a supporting hyperplane for the graph of $u(y)$ but $$y_{n+1}=u(x)+(\nabla u(x)-\delta \nu_x)\cdot (y-x)$$ is not a supporting hyperplane for any $\delta > 0$, where $\nu_x$ denotes the exterior normal to $\p \Omega$ at $x$.
\end{rem}

In general the Lipschitz continuity of the solution can be easily obtained from the boundary data by the use of barriers. Also the quadratic separation assumption \eqref{quad_sep} can be checked in several situations directly from the boundary data, see Proposition 3.2 in \cite{S2}. This is the case for example when the boundary data is more regular i.e \begin{equation}\label{c3}
u|_{\p \Omega}, \p \Omega \in C^3, \quad \mbox{and $\Omega$ is uniformly convex.}
\end{equation}
As a consequence of Theorem \ref{w2p} we obtain
\begin{cor}\label{c1} Let $u:\overline \Omega \to \R$ be a continuous solution of the Monge-Ampere equation \eqref{MA} that satisfies \eqref{f_osc} and \eqref{c3}. Then $$\|u\|_{W^{2,p}(\Omega)} \le C,$$ with $C$ depending on $u|_{\p \Omega}$, $\p \Omega$, $\rho$, $\lambda$, $\Lambda$, $p$ and $n$.
\end{cor}

In particular, if $u$ solves the Monge-Ampere equation \eqref{MA} with $f \in C(\overline \Omega)$ and \eqref{c3} holds, then $u\in W^{2,p}(\Omega)$ for any $p< \infty$.

The assumptions on the boundary behavior of $u$ and $\p \Omega$ in Theorem \ref{w2p} and Corollary \ref{c1} seem to be optimal. Wang in \cite{W} gave examples of solutions $u$ to \eqref{MA} with $f=1$ and either $u\in C^{2,1}$ or $\p \Omega \in C^{2,1}$, that do not belong to $W^{2,p}(\Omega)$ for large values of $p$.

The proof of Theorem \ref{w2p} is based on a localization theorem for the Monge-Ampere equation at boundary points which was proved in [S1], [S2]. It states that under natural local assumptions on the domain and boundary data, the sections
  $$S_h(x_0)=\{x \in \overline \Omega | \quad u(x)<u(x_0)+\nabla u(x_0) \cdot (x-x_0)+h\},$$
with $x_0 \in \p \Omega$ are ``equivalent" to ellipsoids centered at $x_0$. We give its precise statement below.

Assume for simplicity that
\begin{equation}\label{om_ass}
B_\rho(\rho e_n) \subset \, \Omega \, \subset \{x_n \geq 0\} \cap B_{\frac 1\rho},
\end{equation}
for some small $\rho>0$, that is $\Omega \subset (\R^n)^+$ and $\Omega$ contains an interior ball tangent to $\p \Omega$ at $0.$
Let $u : \overline \Omega \rightarrow \R$ be continuous, convex, satisfying
\begin{equation}\label{eq_u}
\det D^2u =f, \quad \quad 0<\lambda \leq f \leq \Lambda \quad \text{in $\Omega$}.
\end{equation}
After subtracting a linear function we also assume that
\begin{equation}\label{eq_u1}
\mbox{$x_{n+1}=0$ is the tangent plane to $u$ at $0$,}
\end{equation}
in the sense that $$u \geq 0, \quad u(0)=0,$$
and any hyperplane $x_{n+1}= \delta x_n$, $\delta>0$, is not a supporting plane for $u$.

Theorem \ref{main_loc} shows that if the boundary data has quadratic growth near $\{x_n=0\}$ then, each section of $u$ at $0$
$$S_h :=S_h(0)= \{x \in \overline \Omega : \quad u(x) < h \},$$
is equivalent to a half-ellipsoid centered at 0.

\begin{thm}[Localization theorem]\label{main_loc} Assume that $\Omega$, $u$ satisfy \eqref{om_ass}-\eqref{eq_u1} above and,
$$\rho |x|^2 \leq u(x) \leq \rho^{-1} |x|^2 \quad \text{on $\p \Omega \cap \{x_n \leq \rho\}.$}$$
Then, for each $h<c_0$ there exists an ellipsoid $E_h$ of volume $\omega_n h^{n/2}$ such that
$$c_0E_h \cap \overline \Omega \, \subset \, S_h \, \subset \, c_0^{-1}E_h \cap \overline \Omega.$$

Moreover, the ellipsoid $E_h$ is obtained from the ball of radius $h^{1/2}$ by a
linear transformation $A_h^{-1}$ (sliding along the $x_n=0$ plane)
$$A_hE_h= h^{1/2}B_1$$
$$A_h(x) = x - \nu x_n, \quad \nu = (\nu_1, \nu_2, \ldots, \nu_{n-1}, 0), $$
with
$$ |\nu| \leq c_0^{-1} |\log h|.$$
The constant $c_0>0$ above depends on $\rho, \lambda, \Lambda$,and  $n$.
\end{thm}

\section{Proof of Theorem \ref{w2p}}

We start by remarking that under the assumptions of Theorem \ref{main_loc} above, we obtain that also the section $\{u<h^{1/2}x_n\}$ has the shape of $E_h$.

Indeed, since $$S_h \subset c_0^{-1} E_h \subset \{x_n \le c_0^{-1} h^{1/2} \} $$ and $u(0)=0$, we can conclude from the convexity of $u$ that the set
$$F:=\{x \in \overline \Omega| \quad u<c_0h^{1/2}x_n \} $$ satisfies for all small $h$
\begin{equation}\label{f_sub}
F \subset S_h \cap \Omega,
\end{equation}
 and $F$ is tangent to $\p \Omega$ at $0$. We show that $F$ is equivalent to $E_h$ by bounding its volume by below.

\begin{lem}\label{l1} We have
$$|F| \ge c|E_h|$$
for some $c>0$ small depending on $\rho$, $\lambda$, $\Lambda$, $n$.
\end{lem}

\begin{proof}
From Theorem \ref{main_loc}, there exists $y \in \partial S_{\theta h}$ such that $y_n \ge c_0(\theta h)^{1/2}$. We evaluate $$v:=u-c_0 h^{1/2}x_n, $$ at $y$ and find $$v(y) \le \theta h - c_0 h^{1/2} c_0 (\theta h)^{1/2} \le -\delta h,$$ for some $\delta>0$ provided that we choose $\theta$ small depending on $c_0$. Since $v=0$ on $\p F$ and $$ \det D^2 v \ge \lambda$$ we have $$|\inf_F v| \le C(\lambda)|F|^{2/n},$$ hence $$ c h^{n/2} \le |F|.$$
\end{proof}

 Next we prove Theorem \ref{w2p}. We denote by $c$, $C$ positive constants that depend on $\rho$, $\lambda$, $\Lambda$, $p$, $n$ and $\|\p \Omega\|_{C^{1,1}}$, $\| u|_{\p \Omega}\|_ { C^{1,1}}$, $\|u\|_{C^{0,1}}$. For simplicity of notation, their values may change from line to line whenever there is no possibility of confusion.

 Given $y\in \Omega$ we denote by
 $$S_h(y):=\{x \in \overline \Omega| \quad u(x) < u(y)+\nabla u(y) \cdot (x-y)\}$$
and let $\bar{h}(y)$ be the maximal value of $h$ such that $S_h(y) \subset \Omega$, i.e
$$\bar{h}(y):=\max\{h\ge 0|\quad S_h(y) \subset \Omega\}.$$

\begin{lem}\label{l2}
Let $y \in \Omega$ and denote for simplicity $\bar h=\bar{h}(y)$. If $\bar h \le c_1$ then
$$S_{\bar h}(y) \subset D_{C\bar h^{1/2}}:=\{x \in \overline \Omega|\quad dist(x,\p \Omega)\le C\bar h^{1/2} \}, $$
and
$$\int_{S_{\bar h/2}(y)}\|D^2 u\|^p \, dx \le C|\log \bar h|^{2p} \, \bar h^{n/2}.$$
\end{lem}

\begin{proof}
Without loss of generality assume $$0 \in \p S_{\bar h}(y) \cap \p \Omega.$$ After subtracting a linear function and relabeling $\rho$ if necessary we may also assume that $u$ satisfies the conditions of Theorem \ref{main_loc} at the origin. This implies
$$S_{\bar h}(y)=\{x\in \overline \Omega| \quad u<\alpha x_n\}, \quad \quad \alpha>0$$ and $\bar h>0$. Since (see \cite{C})
$$c\bar h^{n/2} \le |S_{\bar h}(y)| \le C \bar h^{n/2},$$ we obtain from \eqref{f_sub} and Lemma \ref{l1}
$$c\bar h^{1/2} \le \alpha \le C \bar h^{1/2},$$
and that $S_{\bar h}(y)$ is equivalent to an ellipsoid $\tilde E:=E_{C\bar h}$. This implies that
$$y+c\tilde E \subset S_{\bar h}(y) \subset C\tilde E,$$
with
$$\tilde E =\bar h^{1/2}A^{-1}B_1, \quad Ax:=x-\nu x_n,$$
and  $$ \nu_n=0, |\nu| \le C |\log \bar h|.$$

The inclusion above implies $$S_{\bar h}(y) \subset \{x_n \le C\bar h^{1/2}\} \subset D_{C\bar h^{1/2}}.$$
The rescaling $\tilde u: \tilde \Omega \to \R$ of $u$
$$\tilde u(x):=\frac {1}{ \bar h} u(y+\bar h^{1/2}A^{-1}x),$$
satisfies
$$\tilde \Omega=\bar h^{-1/2}A(\Omega - y), \quad \quad \tilde S_1(0)=\bar h^{-1/2}A(S_{\bar h}(y) - y)\subset \tilde \Omega,$$
and
$$B_c \subset \tilde S_1(0) \subset B_C,$$
where $\tilde S_1(0)$ represents the section of $\tilde u$ at the origin at height 1.
Moreover,
$$\det D^2\tilde u(x)=\tilde f(x):=f(y+\bar h^{1/2}A^{-1}x)$$
and \eqref{f_osc} implies
$$|\log \tilde f(x)-\log \tilde f(z)|\le \eps(p) \quad \quad \forall x,z \in \tilde S_1(0)$$
since $$|\bar h^{1/2}A^{-1}(x-z)| \le C \bar h^{1/2}|\log \bar h| \le \rho.$$
The interior $W^{2,p}$ estimate for $\tilde u$ in $\tilde S_0(1)$ gives
$$\int_{\tilde S_{1/2}(0)}\|D^2 \tilde u\|^p \, dx \le C,$$
hence
$$\int_{S_{\bar h/2}(y)}\|D^2 u\|^p \, dx=\int_{\tilde S_{1/2}(0)}\|A^T \, D^2 \tilde u \, A\|^p \, \bar h^{n/2}\, dx \le C  |\log \bar h|^{2p} \, \bar h^{n/2}.$$

\end{proof}

We also need the following covering lemma.

\begin{lem}[Vitali covering]\label{l3}
There exists a sequence of disjoint sections $$S_{\delta \bar h_i}(y_i), \quad \quad \bar h_i=\bar h(y_i),$$ such that $$\Omega \subset \bigcup_{i=1}^\infty S_{\bar h_i/2}(y_i),$$
where $\delta>0$ is a small constant that depends only on $\lambda$, $\Lambda$ and $n$.
\end{lem}

\begin{proof}
We choose $\delta$ such that if $$S_{\delta \bar h(y)}(y) \cap S_{\delta \bar h(z)}(z) \ne \emptyset \quad \quad \mbox{and} \quad 2\bar h(y) \ge \bar h(z),$$
then $$S_{\delta \bar h(z)}(z)\subset S_{\bar h(y)/2}(y).$$
The existence of $\delta$ follows from the engulfing properties of sections of solutions to Monge-Ampere equation \eqref{MA} with bounded right hand side (see \cite{CG}).

Now the proof is identical to the proof of Vitali's covering lemma for balls. We choose $S_{\delta \bar h_1}(y_1)$ from all the sections $S_{\delta \bar h(y)}(y)$, $y \in \Omega$ such that $$\bar h(y_1) \ge \frac 1 2 \sup_{y}\bar h(y),$$ then choose $S_{\delta \bar h_2}(y_2)$ as above but only from the remaining sections $S_{\delta \bar h(y)}(y)$ that are disjoint from $S_{\delta \bar h_1}(y_1)$, then $S_{\delta \bar h_3}(y_3)$ etc. We easily obtain

$$\Omega=\bigcup_{y\in \Omega} S_{\delta \bar h(y)}(y) \subset \bigcup_{i=1}^\infty S_{\bar h_i/2}(y_i).$$

\end{proof}

{\it End of proof of Theorem \ref{w2p}}

$$\int_\Omega \|D^2 u\|^p \, dx \le \sum_i \int_{S_{\bar h_i/2}(y_i)} \|D^2 u\|^p \, dx.$$

There are a finite number of sections with $\bar h_i \ge c_1$ and, by the interior $W^{2,p}$ estimate, in each such section we have

$$\int_{S_{\bar h_i/2}(y_i)} \|D^2 u\|^p \, dx \le C.$$

Next we consider the family $\mathcal{F}_d $ of sections $S_{\bar h_i/2}(y_i)$ such that $$ d/2 < \bar h_i \le d$$ for some constant $d \le c_1$. By Lemma \ref{l2} in each such section  $$\int_{S_{\bar h_i/2}(y_i)} \|D^2 u\|^p \, dx \le C |\log d|^{2p} |S_{\delta \bar h_i}(y_i)|,$$
and since $$S_{\delta \bar h_i}(y_i) \subset D_{Cd^{1/2}}$$ are disjoint we find
$$\sum_{i \in \mathcal{F}_d} \int_{S_{\bar h_i/2}(y_i)} \|D^2 u\|^p \, dx \le C |\log d|^{2p} d^{1/2}.$$
We add these inequalities for the sequence $d=c_12^{-k}$, $k=0,1,\dots$ and obtain the desired bound.


\begin{thebibliography}{9999}

\bibitem[C]{C} Caffarelli L., Interior $W^{2,p}$ estimates for solutions of Monge-Ampere equation, {\it Ann. of Math.} {\bf 131} (1990), 135-150.

\bibitem[CG]{CG} Caffarelli L., Gutierrez C., Properties of solutions of the linearized Monge-Ampere equations, {\it Amer. J. Math.} {\bf 119} (1997), 423-465.

\bibitem[S1]{S1} Savin O., A localization property at the boundary for Monge-Ampere equation, Preprint arxiv:1010.1745.

\bibitem[S2]{S2} Savin O, Pointwise $C^{2,\alpha}$ estimates at the boundary for Monge-Ampere equation, Preprint arxiv:1101.5436.

\bibitem[W]{W} Wang X.-J., Regularity for Monge-Ampere equation near the boundary, {\it Analysis} {\bf 16} (1996) 101-107.

\end{thebibliography}
\end{document}